\ifx\LocalElevenPtArticleFormatLoaded\undefined
\documentclass[a4paper,11pt]{article} 
\usepackage{theorem} 
\usepackage{amsmath,amscd,amssymb} 
\usepackage{latexsym} 
\usepackage[german,english]{babel} 
\usepackage[applemac]{inputenc}
\usepackage{float}
\usepackage{epsfig}
\usepackage{xypic}
\fi

\newtheorem{theorem}{Theorem}[section]
\usepackage{theorem}

\newtheorem{lemma}[theorem]{Lemma}

\newtheorem{proposition}[theorem]{Proposition}
\newtheorem{corollary}[theorem]{Corollary}
\theorembodyfont{\rm}
\newtheorem{remark}[theorem]{Remark}

\newcommand{\AAA}{{\mathbb{A}}}

\newcommand{\FF}{{\mathbb{F}}}
\newcommand{\PP}{{\mathbb{P}}}
\newcommand{\QQ}{{\mathbb{Q}}}

\newcommand{\ZZ}{{\mathbb{Z}}}

\newenvironment{Proof}{\begin{ProofwCaption}{Proof}}{\end{ProofwCaption}}
\newenvironment{Proof*}[1]{\begin{ProofwCaption}{{#1}}}{\end{ProofwCaption}}
\newenvironment{ProofwCaption}[1]%
  {\addvspace\theorempreskipamount \noindent{\it #1.}\rm}%
  {\qed \par \addvspace\theorempostskipamount}
\newcommand{\qedsymbol}{\mbox{$\Box$}}
\newcommand{\qed}{\quad\qedsymbol}
\DeclareMathOperator{\Frob}{Frob}

\input xy
\xyoption{all}

\begin{document}

\title{On the motive of Kummer varieties associated to $\Gamma_1(7)$ -
Supplement to the paper: The modularity of certain non-rigid Calabi-Yau threefolds
(by R.~Livn\'e and N.~Yui)}
\author{K.~Hulek and H.~Verrill}
\date{}
\maketitle

\begin{abstract}
In their paper \cite{LY} Livn\'e and Yui discuss several examples of non-rigid Calabi-Yau
varieties which admit semi-stable $K3$-fibrations with $6$ singular fibres over a base which is
a rational modular curve. They also establish the modularity of the $L$-function of these
examples. The purpose of this note is to point out that the examples which were listed in
\cite{LY}, but which do not lead to semi-stable fibrations, are still modular in the sense
that their $L$-function is associated to modular forms. We shall treat
the case associated to the group
$\Gamma_1(7)$ in detail, but our technique also works in the other cases given in \cite{LY}.
We shall also make some comments concerning the
Kummer construction for fibre products of elliptic surfaces in general.
\end{abstract}

\section{Introduction}

In their paper \cite{LY} Livn\`e and Yui consider Calabi-Yau varieties which possess a 
non-constant semi-stable 
$K3$-fibration with $6$ singular fibres, which is the minimal number (Arakelov-Yau bound) of 
singular fibres of such a fibration (\cite{STZ}). 
In this case the base curve must be a rational modular curve.
They start with the list of the (up to conjugacy) $9$ possible torsion free genus $0$ congruence 
subgroups of index $24$ of $\operatorname{PSL}(2,\ZZ)$. These $9$ cases separate into two types, 
depending on whether the group
is a subgroup of $\Gamma(2)$ or, equivalently, whether the curve 
of $2$-torsion points decomposes into four sections ($4$ cases) or not ($5$ cases). The first 
situation leads to
semi-stable $K3$ fibrations by performing the following Kummer construction: if $Y$ is the 
universal elliptic curve and $E$ is an elliptic curve 
then $X'=(Y \times E)/ \iota$ (where $\iota$ is the map given by
$x \mapsto -x$) is a singular threefold. The fibre structure of $Y$ induces a fibration on $X'$, 
whose general fibre is a (singular) Kummer surface. 
Blowing up the fixed point set of $\iota$ gives the desired smooth Calabi-Yau variety $X$. 
Livn\'e and Yui also show that these non-rigid Calabi-Yau varieties are modular in the 
sense that their $L$-function is modular. The point of this note is to make some general statements
about Kummer fibrations and to point out that 
the modularity statement is also true in the remaining cases. Here we shall restrict 
ourselves to the $\Gamma_1(7)$ case which is particularly interesting. 
The same method, however, also applies to the other cases. The main difference to the examples treated 
in \cite{LY} is that the curve of non-zero $2$-torsion points is an elliptic curve. This leads to 
extra contributions in the motive of $X$.

\section{The Kummer construction and topological considerations}\label{topology}

Let
$r:Y\rightarrow S$ and $r':Y'\rightarrow S$
be two relatively minimal elliptic fibrations with a section over the same base curve. We assume
these to be semi-stable, i.e. all fibres are of type $I_n, n\geq 0$. The
involutions $\iota_Y$ resp. $\iota_{Y'}$ which, on the general fibre, are
given by $x\mapsto -x$ extend to $Y$, resp. $Y'$. We label the 
components of a $I_n$-fibre cyclically by $e_0, \ldots, e_{n-1}$ and in such a way that $e_0$ 
corresponds to the component which meets the $0$-section. Then $\iota$ acts by $e_i \mapsto e_{-i}$.
Hence, if $n$ is even,
two components are fixed, namely $e_0$ and $e_{n/2}$, and the others are exchanged pairwise, whereas for
$n$ odd only one component is fixed, namely $e_0$. Let $O,B \subset Y$, resp. $O',B'\subset Y'$
be the $0$-section and the closure of the non-zero $2$-torsion points.
These are smooth curves. If $n$ is even, $B$ meets $e_0$ and $e_{n/2}$ transversally 
in $1$, resp. $2$ smooth points,
whereas if $n$ is odd then $B$ meets $e_0$
transversally in one point, as well as the intersection of $e_{(n-1)/2}$ and $e_{(n+1)/2}$. In this case the map $B\rightarrow S$ is
branched.

Let $W=Y\times_{S}Y'$. This variety is singular
and has nodes ($A_1$-sin\-gu\-la\-ri\-ties) exactly at pairs of nodes of the singular fibres.
If, say, $Y'=E \times \mathbb{P}^1$
for an elliptic curve $E$ (and this will be the situation in our example), then $W$ is smooth.
We consider the diagonal involution
$$
\iota=(\iota_Y, \iota_{Y'}): W\rightarrow W.
$$
Let $\tilde{W}$ be the (big) resolution of $W$, i.e., the variety where all nodes are replaced
by quadrics. The involution $\iota$ lifts to $\tilde{W}$ where we shall denote it by the same
letter.
The {\em (singular) Kummer family} associated to the pair $(Y,Y')$  is defined by
$$
X'=\tilde{W}/\iota.
$$
This variety is always singular, even if $W$ is smooth.
The general fibre is a Kummer
surface with $16$ nodes. In the analytic
category one can also consider small resolutions $\hat{W}$ of $W$ where the double points are 
replaced by a $\mathbb{P}^1$. These small resolutions are not necessarily projective and it is possible
that $W$ does not possess a small projective resolution. Moreover, even if a small
projective resolution $\hat{W}$ exists, then it is not clear that $\iota$ lifts to $\hat{W}$.
For a discussion of this issue see \cite{Sch}. In view of the arithmetic applications we have
in mind we are not concerned with the existence of small projective resolutions at this point and
hence we will work with either $W$ if this is smooth or with a big resolution $\tilde{W}$ otherwise.

The singularities of $X'$ come from fixed points
of $\iota$ on $W$. One has to distinguish two cases, namely:
\begin{enumerate}
\item[(1)] There is no fibre $I_n\times I_m$ with both $n,m$ odd.
\item[(2)] There is such a fibre.
\end{enumerate}
The fixed locus of $\iota$ is given by
$$
D := (O+B)\times_S(O'+B')\subset W.
$$
In the case (1)
the curve $D$ does not go through any of the singularities of $W$,
in case (2) it does. In the first case
one obtains a desingularization of $X'$ in a way exactly like the usual
desingularization of Kummer surfaces: Blowing up $\tilde{W}$ along the curve $D$ replaces $D$ by
a $\mathbb{P}^1$-bundle over $D$ whose rulings have
relative degree $-2$. We shall denote the resulting $3$-fold by $Z$. Then
$$
X= Z/ \iota
$$
is a smooth variety, fibred over the base whose general fibre is a smooth
$K3$-surface. Alternatively we could have obtained $X$ by blowing up $X'$ along
its double curve. We call $X$ the {\em smooth Kummer fibration} associated to the
pair $(Y,Y')$.

In the second case the situation is as follows. Locally, near a singular
point we can choose (analytic) coordinates $x,y$ on $Y$ and
$u,v$ on $Y'$ such that the projection onto the base is given by
$(x,y)\mapsto xy$ and $(u,v)\mapsto uv$. Then the fibre product
$Y\times_S Y'$ is locally isomorphic to the subvariety
given by $xy-uv=0$.
We obtain
$\tilde{W}$ by blowing up the $A_1$-singularity given by
$xy-uv=0$, thus inserting a quadric $Q$.
The strict transform of $B$ is a smooth curve meeting
$Q$ in two points. The involution $\iota$ on $Q$ has $4$
fixed points. The first two of these lie
on the strict transform of $B$, the other two are isolated
singularities. In the quotient $X$ the latter two points give
rise to isolated singularities of type $V_{\frac{1}{2}(1,1,1)}$,
i.e. to cones over Veronese surfaces. These can be resolved by a
$\mathbb{P}^2$ with normal bundle $\mathcal{O}_{\mathbb{P}^2}(-2)$.
In particular, these singularities are rational, and in an arithmetic context they will not
contribute in an essential way to the motive of $X$.

We are interested in examples where $X$ is Calabi-Yau. For the examples which we shall consider
we shall now assume
that $r: Y \to {\mathbb P}^1$ is a $K3$-surface and that $Y' = E \times {\mathbb P}^1$, where $E$
is an elliptic curve.
Note that 
$W=Y \times E$ and the curve $D$ are smooth.

We denote the exceptional $\mathbb{P}^1$-bundle over $D$ by $V$.
Moreover, we denote by $n_+$ the
rank of the $(+1)$-eigenspace of $NS(Y)$ with respect
to $\iota_Y$ and by $n_-$ the rank of the $(-1)$-eigenspace.
Let $g(D)$ be the sum of the genera
of the components of $D$ and let $c(D)$ be the number of components of
$D$.

\begin{theorem} \label{theo:cohomology}
Let $r: Y \to {\mathbb P}^1$ be a $K3$-surface and $Y' = E \times {\mathbb P}^1$.
Then the associated Kummer variety $X$ is a Calabi-Yau threefold whose Hodge numbers are as follows:
\begin{itemize}
\item[{\rm(i)}] $h^{00}(X)=1, \; h^{10}(X)=h^{01}(X)=0$,
\item[{\rm(ii)}]$h^{20}(X)=h^{02}(X)=0, \; h^{11}(X)=\rho(X)=
n_++1+c(D)$,
\item[{\rm(iii)}]$h^{30}(X)=h^{03}(X)=1, \; h^{12}(X)=h^{21}(X)=
1+n_-+g(D)$.
\end{itemize}
Moreover, for the Euler numbers
$$
e(X)=\frac{3}{2}e(D)
$$
and
$$
n_+-n_-=\frac{3}{4}e(D)-c(D)+g(D)=\frac{1}{4}e(D).
$$
\end{theorem}
\begin{Proof}
In this case $W$ is smooth and there exists a (up to a scalar) unique
$3$-form on $W$ which descends to $X$. As in the case of
Kummer surfaces one shows that this form has no zeroes on $X$ and
hence $\omega_X= \mathcal{O}_X$ and thus $h^{30}(X)=h^{03}(X)=1$. Since there is no invariant $1$-form on $W$
we can also conclude that $h^{10}(X)=h^{01}(X)=0$ and hence by
Serre duality $h^{20}(X)=h^{02}(X)=0$. This shows that $X$ is a Calabi-Yau
threefold. One immediately obtains from the Hodge diamond that
$$
e(X)=2 (h^{11}(X)-h^{12}(X)).
$$
In order to compute the Betti numbers of $X$
we first note that it is easy to determine the Betti numbers of $W$ from those of $Y$ and
$E$ via the K\"unneth formula.
Let $\pi:Z\rightarrow W$ be the blow-up along $D$. Then
$$
H^\ast(Z) \cong \pi^{\ast}H^{\ast}(W)\oplus H^{\ast}(V)/
({\pi}_{|V})^{\ast} (H^{\ast}(D)).
$$
Moreover, recall that $H^{\ast}(V)$ is generated by the tautological class as a ring over $H^{\ast}(D)$.
{}From this one obtains
$$
h^0(Z)=h^0(W), \; h^1(Z)=h^1(W)
$$
$$
h^2(Z)=h^2(W)+ c(D), \; h^3(Z)=h^3(W)+h^1(D).
$$
For the Euler characteristic of $X$ we thus find
\begin{eqnarray*}
e(X)&=&\frac{1}{2}(e(W)-e(D))+e(V)\\
&=& \frac{1}{2}(e(W)-e(D))+2e(D)\\
&=&\frac{1}{2}(e(W)+3e(D))\\
&=&\frac{3}{2}e(D)
\end{eqnarray*}
where the last equality follows from $e(W)=0$.

In order to compute $h^2(X)=\rho (X)$ we have to determine
the $\iota$-invariant divisors on $Z$. The divisors on
$Y$ are spanned by a general fibre $F$, the $0$-section $O$,
and the components $e_i^j$ of the singular fibres which do not
meet the $0$-section. The divisors on $W$ are spanned by
taking the product of these divisors times $E$ and by
$Y\times\left\lbrace 0\right\rbrace$. The pullback of this divisor under
$\pi$ and the components of the exceptional divisors are
clearly invariant under $\iota$. This accounts for the summand
$1 + c(D)$ in the formula for $h^{11}(X) = \rho (X)$.
If $l\in \operatorname{NS}(Y)^{\pm}$ then $l\times E$ is invariant (anti-invariant)
under $\iota$ and altogether this shows the formula for $h^{11}(X)$.
It remains to determine the invariant part of $H^3(Z)$. We shall first treat the
contribution from $W= Y \times E$.
Let $T \subset H^2(Y)$ be the rank $2$ subspace spanned by
the transcendental cycles. Then $T\otimes H^1(E)$ is $4$-dimensional
and $\iota$ acts by $-1$ on both $T$ and $H^1(E)$. This is clear
for $H^1(E)$. If it did not act by $-1$ on $T$ then $h^{20}(X)
\neq 0$ resp. there would be no $3$-form on $X$. Hence $\iota$
acts by $+1$ on $T\otimes H^1(E)$ and this contributes to
$H^3(X)$. The same argument applies to $\operatorname{NS}(X)^-\otimes H^1(E)$.
Moreover $H^3(V)$ is invariant under $\iota$ and in total we
find that
$$
h^3(X)=4+2n_-+2g(D).
$$
Since we already know that $h^{30}(X)=h^{03}(X)=1$ this
proves the claim for $h^{12}(X)$. The final assertion follows from subtracting (iii) from
(ii).\hfill
\end{Proof}

\begin{remark}
\label{rem:exact_sequence}
The middle cohomology of $X$ fits into an exact sequence
$$
0\rightarrow T\otimes H^1(E) \rightarrow H^3(X,\mathbb{C})
\rightarrow (\operatorname{NS}(Y)^- \otimes H^1(E)) \oplus H^3(V)
\rightarrow 0.
$$
\end{remark}

\begin{proposition} \label{prop:valueofn}
Let $r: Y \to S$ be an elliptic fibration with semi-stable fibres.
The numbers $n_+$ and $n_-$ are given by
$$
n_+ - n_- = 2 + \# \{singular~fibres~I_n ~with~n>1~even \}
$$
and $n_+ + n_- = \rho(Y)$.
\end{proposition}
\begin{Proof}
The N\`eron-Severi group $\operatorname{NS}(Y)$ is generated by a general fibre $F$, the $0$-section $O$
and the components $e^j_i$ of the singular fibres $I_n$ where $j$ runs through the cusps and
$i= 0, \ldots, n-1$. 
Then $\iota$ acts on these components by $e^j_i \mapsto e^j_{-i}$. This means the following: if
$n$ is odd then this fibre contributes equally to $n_+$ and $n_-$, whereas for even $n$ we have
one more invariant component (namely $e^j_{n/2}$), than odd components. Since $F$ and the $0$-section $O$
are clearly invariant we obtain the first claim. The second claim follows since  $n_+ + n_- = 2 + \rho(Y)$.
\hfill
\end{Proof}

\section{The example associated to $\Gamma_1(7)$}

From now on we shall concentrate on the example associated to the group $\Gamma_1(7)$. We note, however,
that this method can also be applied to the other examples in \cite{LY} which do not lead to semi-stable
$K3$-fibrations.

Let $Y=S(\Gamma_1(7))$. Then $Y$ has $6$ singular fibres, namely $3$ of type $I_1$ and $3$
of type $I_7$.
By results of Igusa, one knows that the moduli problem for
$\Gamma_1(N)$ for $N>3$ can be represented by a smooth scheme
$Y_1(N)$ over $\ZZ[1/N]$ 
(see \cite[Introduction]{KM} and \cite[Table~(10.9.6), p. 308]{KM}).
Since we shall need this, we shall briefly sketch the construction of a smooth
model defined over $\QQ$.
Tate \cite[p.~195]{Tate}
(and also \cite[case~15., Table~3, p.~217]{kubert})
gives the following equation for the universal elliptic curve
with a group of sections of order $7$
\begin{equation}
y^2+(1+t-t^2)xy+(t^2-t^3)y=x^3+(t^2-t^3)x^2
\label{weier7}
\end{equation}
Note that after a change of variables
$y\mapsto y-\frac{1}{2}(1+t-t^2)x-\frac{1}{2}(t^2-t^3)$ this has the
Weierstrass form
\begin{equation}
4y^2=4x^3+(t^4-6t^3+3t^2+2t+1)x^2+2t^2(t^3-2t^2+1)x+t^4(t-1)^2
\label{weier7B}
\end{equation}
This shows that the curve $B$ of non-zero $2$-torsion points is given by the following 
equation (which was first pointed out to us by N.~Yui) 
\begin{equation}
4x^3+(t^4-6t^3+3t^2+2t+1)x^2+2t^2(t^3-2t^2+1)x+t^4(t-1)^2=0
\label{eqn_for_B}
\end{equation}

\begin{lemma}
\label{lem:L_series_of_B}
The curve $B$ of non-zero $2$-torsion points of $\;Y$ is an 
elliptic curve defined over $\QQ$, with conductor $14$.  The
Mellin transform of the L-series of this curve is the unique
weight $2$ level $14$ modular newform, $\eta(\tau)\eta(2\tau)\eta(7\tau)
\eta(14\tau)$.
\end{lemma}
\begin{Proof}
The curve $B$ parametrizes elliptic curves with a point of order $7$
and a point of order $2$.  This is an irreducible family, and by the
Hurwitz formula it has genus $1$.
In particular, $B$ is isomorphic to the modular curve
$X_1(2,7)=\Gamma_1(2)\cap\Gamma_1(7)\setminus\mathfrak h^*$
(where $\mathfrak h^*$ is the upperhalf complex plane union the cusps).
There is a degree $3$ natural map from this curve to $X_0(2,7)
=\Gamma_0(2)\cap\Gamma_0(7)\setminus\mathfrak h^*$ which parametrizes elliptic
curves together with a subgroup of order $2$ and a subgroup of order $7$.
This curve is isomorphic to $X_0(14)$, which is an elliptic curve.
Hence the map from $B$ to $X_0(14)$ must be an isogeny of
degree $3$.
The fact that the map from $X_1(2,7)$ to $X_0(14)$ is defined over
$\QQ$ follows, viewing $X_1(2,7)$ as a fibre product, from
the fact that the maps $X_1(7)\rightarrow X_0(7)$
(\cite[(4.24)]{E}) and $X_0(2)\rightarrow X_0(1)$, in the
following diagram, are defined over $\QQ$.
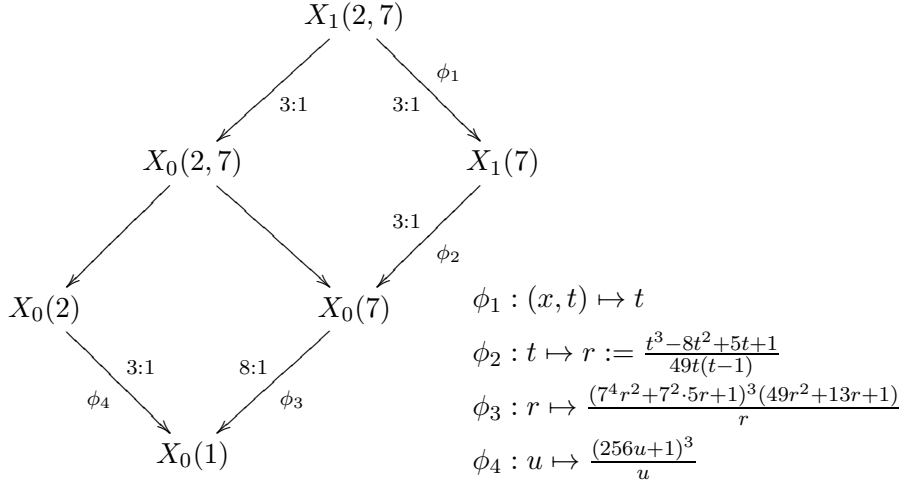
$$
\begin{array}{ll}
\xymatrix@R=1.35cm@C=0.6cm{
&&X_1(2,7)\ar[dr]_{3:1}^{\phi_1}\ar[dl]^{3:1}\\
&X_0(2,7)\ar[dr]\ar[dl]&&X_1(7)\ar[dl]_{3:1}^{
\phi_2}\\
X_0(2)\ar[dr]^{3:1}_{\phi_4}
&&X_0(7)\ar[dl]_{8:1}^{\phi_3}\\
&X_0(1)
}
&
\hspace{-0.6in}
\raisebox{-4.8cm}{$
\renewcommand{\arraystretch}{1.5}
\begin{array}{l}
\phi_1:(x,t)\mapsto t\\
\phi_2: t   \mapsto r:=\frac{t^3-8t^2+5t+1}{49t(t-1)}\\
\phi_3: r   \mapsto \frac{(7^4r^2+7^2\cdot5r+1)^3(49r^2 + 13r +1)}{r}\\
\phi_4: u   \mapsto \frac{(256u+1)^3}{u}
\end{array}$
}
\end{array}
$$
In the diagram $x$ and $t$ are modular functions satisfying the
relation (\ref{eqn_for_B}), and
$u=(\eta(2\tau)/\eta(\tau))^{24}$ and
$r=(\eta(7\tau)/\eta(\tau))^4$ are hauptmodules for
$\Gamma_0(2)$ and $\Gamma_0(7)$ respectively,
taken from
\cite[Table 3]{CN}.  The maps given in the diagram may be found
by comparison of $\QQ$-expansions, and the fact that the degree of
the rational functions equals the index of the corresponding groups.
This implies that $B$ has the same L-series as
$X_0(14)$, which is the given modular form (listed in \cite{martin}).
\hfill
\end{Proof}

As a surface (\ref{weier7})
is isomorphic to $S(\Gamma_1(7))$ away from the singular
fibres.  This example is discussed in more detail by Elkies
\cite[\S4.2]{E}, where the parameter $t$ is given explicitly as a
Hauptmodul for $\Gamma_1(7)$.
The $j$-invariant of this elliptic curve (for both
(\ref{weier7}) and (\ref{weier7B}))
is
$$
\frac
{(t^2 - t + 1)^3(t^6 - 11t^5 + 30t^4 - 15t^3 - 10t^2 + 5t + 1)^3}
{t^7(t-1)^7(t^3 - 8t^2 + 5t +1)},
$$
which implies that, after resolution of singularities,
this
model has $I_7$ fibres at $t=0,1,\infty$, and $I_1$
fibres at the roots of $t^3 - 8t^2 + 5t +1=0$.

The surface (\ref{weier7}) is singular, but it is important for
us to have a smooth resolution defined over $\QQ$, so
we now give an outline of an explicit desingularization
of (\ref{weier7}).
First, we make the change of variables
\begin{equation}
x = x't^2,
y = y't^2(t-1),
z = x'/(t-1) + y' + z',
\label{eqn:changeofvariables}
\end{equation}
which is invertible on the fibres for $t\not=0,1,\infty$.
This gives us
\begin{equation}
\label{eqn_with_new_coords}
 t(t-1)x(x-y)(y+z)+ (t-1)(x-y-z)yz + t(x-y)x z = 0.
\end{equation}
(where we write $x,y,z$
instead of $x',y',z'$).  
This surface (\ref{eqn_with_new_coords}) has $8$ singular
points, namely
$P_1=(1:0:1), P_2=(1:1:0), P_3=(1:0:0)$ on the fibre $t=0$,
$Q_1=(0:0:1), Q_2=(0:1:0), Q_3=(1:1:0)$ over $t=1$ and
$R_1=(0:0:1)$ and $R_2=(0:1:-1)$ over $t=\infty$.
Note that these points also lie on the closure of the sections given by the
points of order $7$.
These points may be resolved by a sequence of blow ups in the singular
points. The first blow up replaces $P_1, P_2, Q_1, Q_2$ and $R_1$ by one line,
and the other points by $2$ lines. The resulting surface still has one singular point, namely
a point on the singular fibre over $t=\infty$, which is infinitesimally close to $R_2$.
Blowing up once more in this point gives rise to a further line and the resulting surface is then
smooth. It has three singular fibres of type $I_1$ and three of type $I_7$ and is isomorphic to
$S(\Gamma_1(7))$. 

The locus which is blown up is defined over $\QQ$, and so the
resulting surface is also defined over $\QQ$
(see the discussion in the proof of Proposition~\ref{prop:normalbundle}).
Moreover, going through the above sequence of blow ups shows that all components of the 
singular fibres are defined over $\QQ$.
This gives us the following result.
\begin{proposition}
\label{prop:NSYgenerators_over_Q}
The group $\operatorname{NS}(Y)$ can be generated by classes of curves
defined over ${\mathbb Q}$.
\end{proposition}
\begin{Proof}
As already mentioned (in the proof of Proposition~\ref{prop:valueofn}),
$\operatorname{NS}(Y)$ is generated by a general fibre, the zero section,
and components of the singular fibres not meeting this section. 
Since all components of the singular fibres and the zero-section are defined
over $\QQ$ this gives the result.
\hfill
\end{Proof}

\begin{proposition}
\label{prop:L_series_of_Y}
The Mellin transform of the
L-series of the summand of the middle cohomology of
$Y$
corresponding to the transcendental lattice
is given by $(\eta(\tau)\eta(7\tau))^3$,
the unique normalized
weight $3$ modular form for $\Gamma_1(7)$,
where
$\eta(\tau)=\prod_{n>0}(1-\exp(2\pi in\tau))$ is the Dedekind $\eta$-function.
\end{proposition}
\begin{Proof}[Sketch]
We first note that modularity as such can be deduced from Livn\'e's 
result \cite[Example 1.6]{livne2}.
One may also apply Serre and Deligne's \cite{D} methods,
(details given by Conrad \cite{C}), as is done
in a related situation in \cite{SY}.
Using the formula for the dimensions of spaces of cusp forms
and applying the transformation properties of
the Dedekind $\eta$-function, one shows that the corresponding  form
is as given.

Alternatively, in order to determine the modular form, we can use the Lefschetz trace formula to compute
the values of the trace of $\Frob_p$ on $H^2(Y,\QQ_\ell)$
by counting points on $Y$
(using e.g., {\sc Magma} \cite{magma}).
These may be
compared with coefficients of $(\eta(\tau)\eta(7\tau))^3$. 
We want to use Livn\'e's result \cite{livne} to show that the
two representations coincide up to semi-simplification.
In order to apply this result we first observe that the determinant is $ \chi_{-7} \chi_{\ell}^2$
(again cf. \cite[Example 1.6]{livne2} and \cite[(2.1)]{LY}). Here $\chi_{\ell}$ is 
the cyclotomic character and $\chi_{-7}$ is the quadratic
character associated to the field $\QQ(\sqrt{-7})$.
We also need that the number of points on $Y$ is even.
This depends on the curve $B$ having an even number of points. Since $B$ and $X_0(14)$ are isogenous 
over $\QQ$ by Lemma \ref{lem:L_series_of_B} it is enough to consider the latter curve.
Looking at a Weierstrass equation one observes that this curve has a non-trivial
$2$-torsion point over $\QQ$, and this is enough to conclude that the number of points over
$\FF_p$ is even for almost all $p$.
This approach
also essentially uses Serre and Deligne's work for the
existence of a modular Galois representation corresponding to
$(\eta(\tau)\eta(7\tau))^3$. \hfill
\end{Proof}

\begin{remark}
The form $g_3$ is related to the unique newform $h_2$ of weight $2$ and level $49$ as follows:  
in terms of representations (at least up to semi-simplification)
$\operatorname{Sym}^2(h_2)=g_3 \otimes \chi_{-7}\chi_{\ell}$.
\end{remark}

For what follows we need to understand the exceptional locus of the blow up
along the fixed point set of $\iota$.

\begin{proposition}\label{prop:normalbundle}
The exceptional locus $V \subset Z$ has the following properties:
\begin{itemize}
\item[\rm{(i)}] $V$ consists of $8$ components. Of these $4$ are
Hirzebruch surfaces $\Sigma_2$ and $4$ are isomorphic to $B \times \mathbb{P}^1$.
\item[\rm{(ii)}] If the $2$-torsion points of $E$ are defined over $\mathbb Q$, then so
are the components of $V$. In this case the isomorphism of the $4$ non-rational components with
$B \times \mathbb{P}^1$ is defined over $\mathbb Q$.
\end{itemize}
\end{proposition}
\begin{Proof}
We first note that the fix curve $D$ of $\iota$ decomposes as
$$
D=(O+B)\times_{\mathbb{P}^1}(O_1+O_2+O_3+O_4)
$$
where the $O_i$ correspond to the four $2$-torsion points of $E$,
and $B$ is an elliptic curve by Lemma~\ref{lem:L_series_of_B}.

We start with the geometric statements.
Let $D'$ be one of the components of $D$.
Then the component $V'$ of $V$ which lies over $D'$ is
isomorphic to the projectivized normal bundle $\mathbb{P}(N_{D'/W})=\mathbb{P}(N_{D'/Z})$. 
Since $Y$ is a $K3$-surface, and hence has trivial
canonical bundle, the adjunction formula shows that $N_{D'/Y} \cong \mathcal{O}_{\mathbb{P}^1}(-2)$
if $D'$ is rational and $N_{D'/Y} \cong \mathcal{O}_{D'}$ if $D'$ is elliptic. Hence
$N_{D'/W} \cong \mathcal{O}_{\mathbb{P}^1}(-2) \oplus \mathcal{O}_{\mathbb{P}^1}$ or
$N_{D'/W} \cong 2\mathcal{O}_{D'}$.

Now assume that the $2$-torsion points of $E$ are defined over $\mathbb Q$. Then the same is true for
all components of the curve $D$. The fact that the components of the exceptional divisor are then
defined over $\mathbb Q$, follows from the general fact that blowing up is compatible with base
change (cf. the proof of \cite[Theorem 1.1.9]{Liu}). In the case of the non-rational components we saw from the
adjunction formula that the normal bundle is trivial, i.e. we
have isomorphisms $N_{D'({\mathbb C})} \cong 2{\mathcal{O}}_{D'({\mathbb C})}$ and, in particular,
$H^0({{\cal{H}}om}_{{\mathcal{O}}_{D'({\mathbb C})}}(N_{D'({\mathbb C})}, {\mathcal{O}}_{D'({\mathbb C})}))$ is
$2$-dimensional. Since cohomology commutes with the flat base extension ${\mathbb Q} \hookrightarrow
{\mathbb C}$ (\cite[Proposition III.9.3]{Ha}) this implies that the ${\mathbb Q}$ vector space
$H^0({{\cal{H}}om}_{{\mathcal{O}}_{D'({\mathbb Q})}}(N_{D'({\mathbb Q})}, {\mathcal{O}}_{D'({\mathbb Q})}))$
is also $2$-dimensional and in particular there is a homomorphism
$N_{D'({\mathbb Q})} \to 2{\mathcal{O}}_{D'({\mathbb Q})}$ which, after tensoring with ${\mathbb C}$, is an
isomorphism. Using once more flatness of ${\mathbb Q} \hookrightarrow {\mathbb C}$ it follows that
$N_{D'({\mathbb Q})} \cong 2{\mathcal{O}}_{D'({\mathbb Q})}$.
\hfill
\end{Proof}




\begin{theorem}\label{invariantcohomology}
Assume that the elliptic curve $E$, as well as its $2$-torsion points, are defined over
${\mathbb Q}$. Then the $\iota$-invariant part of the middle cohomology of the blow up $Z$ of
$Y(\Gamma_1(7)) \times E$ along the fixed locus of $\iota$ is modular.
More precicely
$$
L(H^3_{\scriptsize{\mbox{\'et}}}(Z)^{\iota},s)\circeq L(g_3 \otimes g^E_2,s)L(g^E_2,s-1)^9
L(g^B_2,s-1)^4
$$
where $g_3$ is the weight $3$ form associated to the
transcendental lattice $T$ and $g^E_2$ and $g^B_2$ are the weight $2$ cusp forms associated to the elliptic
curves $E$ and $B$ respectively.
In terms of the Dedekind $\eta$-function, we have
\begin{eqnarray*}
g_3(\tau) &=& (\eta(\tau)\eta(7\tau))^3=
 q - 3q^2 + 5q^4 - 7q^7 - 3q^8 + 9q^9 - 6q^{11}+\cdots,\\
g_2^B(\tau) &=& \eta(\tau)\eta(2\tau)\eta(7\tau)\eta(14\tau)
=q - q^2 - 2q^3 + q^4 + 2q^6 + q^7 - q^8 + \cdots.
\end{eqnarray*}
(Here $\circeq$ denotes equality of the $L$-series up to finitely many primes,
and $q=\exp(2\pi i\tau)$.)
\end{theorem}
\begin{Proof}
We first note that 
$n_+(Y)=11$ and $n_{-}(Y)=9$.
This follows from Proposition~\ref{prop:valueofn} and the fact that
there are $3$ fibres of type $I_1$ and $3$
of type $I_7$. 
Now recall from Remark~\ref{rem:exact_sequence},
that we have the following exact sequence of cohomology groups
$$
0\rightarrow T\otimes H^1(E) \rightarrow H^3(Z,\mathbb{C})^{\iota}
\rightarrow (\operatorname{NS}(Y)^- \otimes H^1(E)) \oplus H^3(V)
\rightarrow 0
$$
where
$$
H^3(V) \cong H^3(E \times {\mathbb P}^1)^4.
$$
By Proposition~\ref{prop:NSYgenerators_over_Q}
the standard
generators of $\operatorname{NS}(Y)$ are defined over ${\mathbb Q}$.
Moreover, by Proposition \ref{prop:normalbundle} all components of the
exceptional locus $V$ are defined over $\mathbb Q$ and there are $4$ components with non-vanishing third
cohomology, each isomorphic to $E \times {\mathbb P}^1$ (over $\mathbb Q$).
Hence the above sequence can be read as a sequence of \'etale cohomology groups and this implies the
assertion.
We saw that $g_2^B$ and $g_3$ have the given $\eta$-product expressions
in Proposition~\ref{prop:L_series_of_Y} and Lemma~\ref{lem:L_series_of_B}.
\hfill
\end{Proof}

As an immediate application we obtain the $L$-series of the middle 
cohomology of the Kummer variety $X$. 



\begin{corollary}
\label{cor:L-series_on_X}
Assume that the elliptic curve $E$ as well as its $2$-torsion points are defined over
${\mathbb Q}$. Then the Kummer variety $X=X(\Gamma_1(7))$ is modular,
in the sense that its $L$-series
of its middle cohomology
is given in terms of modular forms as follows, where $g^E_2, g_2^B$ and
$g_3$ are as in
Theorem~\ref{invariantcohomology}.
$$
L(H^3_{\scriptsize{\mbox{\'et}}}(X),s) \circeq 
L(H^3_{\scriptsize{\mbox{\'et}}}(Z)^{\iota},s)\circeq L(g_3 \otimes g^E_2,s)L(g^E_2,s-1)^9
L(g^B_2,s-1)^4.
$$
\end{corollary}
\begin{Proof}
Since $Z$ and the involution $\iota$ are defined over ${\mathbb Q}$, general theory
implies that $X$ is also defined over ${\mathbb Q}$ and
$L(H^3_{\scriptsize{\mbox{\'et}}}(X),s)
\circeq L(H^3_{\scriptsize{\mbox{\'et}}}(Z)^{\iota},s)$.
\hfill
\end{Proof}

\begin{remark}
It follows from the work of Kim and Shahidi \cite{KS} that 
$L(g_3 \otimes g^E_2,s)$ has the expected analytic properties.
\end{remark}

\section{Point counting}\label{pointcounting}

We will now make a verification of the
above result (Corollary~\ref{cor:L-series_on_X}) about
the $L$-series of the Kummer variety $X$ by a counting argument.
Suppose that we have $q$-expansions
$$g_2^B=\sum_{n\ge1} b_nq^n,\>
g_2^E=\sum_{n\ge1} c_nq^n,\> \text{ and }
g_3=\sum_{n\ge1} a_nq^n,$$ where $g_2^B$, $g_2^E$ and $g_3$ are as in
Theorem~\ref{invariantcohomology}.
The following table gives us the traces on the cohomology of $X$
predicted by Corollary~\ref{cor:L-series_on_X},
Proposition~\ref{prop:NSYgenerators_over_Q},
and the fact that
$\iota$ is defined over $\QQ$.
$$
\begin{array}{ccc}
&\text{dimension} &\text{trace of }\Frob_p\\
H_6 & 1 & p^3\\
H_5 & 0 & 0\\
H_4 & 20& 20p^2\\
H_3 & 30& a_pc_p + 9pc_p + 4pb_b\\
H_2 & 20& 20p\\
H_1 & 0 & 0\\
H_0 & 1 & 1.
\end{array}
$$
{}From the Lefschetz fixed trace formula this tells us that
\begin{equation}
\label{eqn:predicted_point_count}
n_p=p^3 + 20p^2 - (a_pc_p + 9pc_p + 4pb_b) + 20p + 1.
\end{equation}
We will now compute the number of points in another way and
show that we still get the same result.

In the following
our notation remains as in section \ref{topology}.
We will suppose that the characteristic is not $2$ or $3$,
and that all two-torsion points of $E$ are defined over $\QQ$.
Suppose that $Y$ and $E$ are given by the following
Weierstrass equations
in affine space,
\begin{eqnarray*}
Y: y^2=p_1(x,t), &&
E: y_2^2 = p_2(x_2),
\end{eqnarray*}
where $$p_1(x,t):=
x^3+ (t^4-6t^3+3t^2+2t+1)x^2  + 8t^2(t^3-2t^2+1)x + 16t^4(t-1)^2$$
is obtained from (\ref{weier7B}) by scaling $x$ and $y$.
The involution on the product is now given by
$$\iota:(x,y,x_2,y_2,t)\mapsto(x,-y,x_2,-y_2,t).$$
Now define a rational map $Y\times E\rightarrow\AAA^5$ by
\begin{eqnarray*}
\phi:(x,y,x_2,y_2,t)&\mapsto&(x,y^2,x_2,yy_2,t).
\end{eqnarray*}
Then the points in the image of this map lie on a variety $X'$ having
equation
\begin{eqnarray}
\label{eqn:equationsforX'}
X': & w=p_1(x,t), &
w_2^2=p_1(x,t)p_2(x_2),
\end{eqnarray}
where $\AAA^5$ is taken to have coordinates $x,w,x_2,w_2,t$.
Since $\phi$ is surjective onto $X'$, and identifies exactly points
identified by $\iota$, $X'$ gives a model of an affine piece of
the quotient, away from the
$I_7\times E$ fibres.  The situation at the $I_7\times E$ fibres is
as follows. Let $e_0, e_1, \ldots , e_6$ be the components of the $I_7$ fibres
enumerated as usual in such a way that $\iota$ acts by $e_i \mapsto e_{-i}$.
Then the products $e_i \times \PP^1$ for $i=1,6$ and $i=2,5$ are identified
and give rise to two copies of $E \times \PP^1$ in
the quotient. The component $e_0 \times E$ is mapped to itself and
$\iota$ has $8$ fixed points on this surface. The resulting quotient is a rational surface with
$8$ nodes which are resolved by blowing up the fixed locus of $\iota$. Finally,
the two surfaces $e_3 \times E$ and $e_4 \times E$ are interchanged. Since,
however, the intersection point of $e_3$ and $e_4$ is fixed it follows that the
resulting quotient is a non-normal surface which is singular along a $\PP^1$.

We can now count the number of points
over $\FF_p$
on the Kummer variety $X$ and find the following expression
\begin{equation}
\label{eqn:summation}
n_p= \#X' + \#X_{\infty} + 6\#A + 3\#B + \#C + 2\#F + \#V - \#D 
\end{equation}
where the terms of this equation have the following meaning.
First of all $\#X'$ are all the points on the affine model counted by the 
Legendre symbol, i.e.
\begin{equation}
\label{eqn: pointsonX'}
\#X' = \sum_{x,x_2,t\in\FF_p}
\left(\left(\!\!\!
\frac{\>\;p_1(x,t)p_2(x_2)\>\;}{p}
\!\!\!\right)
+1
\right).
\end{equation}
$A$ is the surface $\AAA^1 \times E$. This term counts the contribution
from the ruled surfaces coming from the components $e_i \times E, i=1,2$ from the 
three $I_7$ fibres. $B$ comes from the component $e_3 \times E$ and $C$ comes from
$e_0 \times E$. The terms for $A$ and $B$ are counted three times, corresponding
to the three fibres over $t=0,1,\infty$. The term $C$ refers to $t= \infty$ only, since
the corresponding terms for $t=0,1$ are already taken care of by the term $X'$ (up to a
correction given by $F$, see below).
Note that we have taken care to count the fibres over the intersection
points $e_i \cap e_{i+1}$ only once. Recall that
\begin{equation}
\label{eqn:numberonE}
\#E(\FF_p)=
\sum_{x_2\in\FF_p}
\left(\left(\!\!\!
\frac{\>\;p_2(x_2)\>\;}{p}
\!\!\!\right)
+1
\right)
+1
=
p+1-c_p.
\end{equation}
We now find from our geometric discussion and the fact that everything is defined
over $\QQ$ that
\begin{equation}
\label{eqn:A}
\#A=\#(\AAA^1 \times E)=p(p-c_p+1) 
\end{equation}
and 
\begin{equation}
\label{eqn:B}
\#B=\#(\AAA^1 \times E) - \#E + \#\PP^1 =p(p-c_p+1)+c_p. 
\end{equation}
The surface which is the blow up of the quotient of $e_0 \times E$
is a smooth rational surface with Picard number $10$. Moreover the 
N\`eron-Severi group is generated by elements defined over $\QQ$. Hence
it has $1 + 10p + p^2$ points. It contains $8$ lines which belong to the
exceptional divisor $V$. Removing these, but counting the points on $D$, we get 
the contribution
\begin{equation}
\label{eqn:C}
\#C=p^2 + 2p + 1. 
\end{equation}
This is counted once, namely for the $I_7$-fibre over $\infty$.
The term $F$ comes from the following effect. The surfaces  coming 
from $e_0 \times E$ over the $t=0,1$ are counted in the expression
given by the Legendre symbol. However, this needs a correction. In 
this affine model the elliptic fibre is a nodal cubic and over the
node we have a $\PP^1$. However, when we blow up this node gets resolved and we 
have two fibres isomorphic to $E$ which are identified by the involution. This
means that we need the correction term
\begin{equation}
\label{eqn:F}
\#F=\#E - \# \PP^1= -c_p. 
\end{equation}
For the exceptional locus we find
\begin{equation}
\label{eqn:V}
\#V - \#D = 4p(\#B + \#\PP^1) = 4p(p-b_p + 1 + p +1). 
\end{equation}
This can be seen as follows: the ruled surfaces over the elliptic components of $D$ are
isomorphic (over $\QQ$) to $B \times \PP^1$ and the other components are Hirzebruch surfaces
$\Sigma_2$ which have a basis of the N\`eron-Severi group which is defined over $\QQ$.
Finally $X_{\infty}$ means points
where $t\not=\infty$, but $x$ or $x_2=\infty$. We find that  
\begin{equation}
\label{eqn:pointsatinfinity}
\#X_\infty= \#(Y_{t \neq \infty} \cup (E\times\AAA^1))/\iota
=2\#(\PP^1\times\AAA^1) - \#\AAA^1 = 2p^2 + p.
\end{equation}
We can now rewrite the number of points as
\begin{eqnarray}
\nonumber
n_p&=&\sum_{x,x_2,t\in\FF_p}
\left[
\left(\left(\!\!\!
\frac{\>\;p_1\>\;}{p}
\!\!\!\right)
+1
\right)
\left(\left(\!\!\!
\frac{\>\;p_2\>\;}{p}
\!\!\!\right)
+1
\right)
-
\left(\!\!\!
\frac{\>\;p_1\>\;}{p}
\!\!\!\right)
-
\left(\!\!\!
\frac{\>\;p_2\>\;}{p}
\!\!\!\right)
\right] +
\\
&&
+ \#X_{\infty} + 6\#A + 3\#B + \#C + 2\#F + \#V - \#D 
,
\label{eqn:pointcount}
\end{eqnarray}
where $p_1=p_1(x,t)$ and $p_2=p_2(x_2)$.
Finally note that from the Lefschetz fixed point theorem,
bearing in mind that $\rho(Y)=20$, and that $\operatorname{NS}(Y)$ is generated by
classes of curves defined over $\QQ$
(Proposition~\ref{prop:NSYgenerators_over_Q}), we have
\begin{equation}
\label{eqn:numberonY}
\#Y(\FF_p)=
\sum_{x,t\in\FF_p}
\left(\left(\!\!\!
\frac{\>\;p_1(x,t)\>\;}{p}
\!\!\!\right)
+1
\right)
+20p
=
p^2 + 20p + a_p + 1.
\end{equation}
Note that in these expressions we have taken into account the
points at infinity on $E$ and on the fibres of $Y$,
and the components of the $I_7$ fibres of
$Y$, which are not counted by the factor involving the Legendre symbol.
Using (\ref{eqn:numberonE}) and  (\ref{eqn:numberonY}) together 
with (\ref{eqn:pointcount}), we obtain
\begin{eqnarray*}
\nonumber
n_p&=&
(\#E(\FF_p)-1)(\#Y(\FF_p)-20p)\\
&&
-
(\#E(\FF_p)-1-p)p^2
-
(\#Y(\FF_p)-20p-p^2)p
\\
&&
+\#X_{\infty} + 6\#A + 3\#B + \#C + 2\#F + \#V - \#D 
\\
&=&(p-c_p)(p^2 + a_p + 1)
-(p-c_p -p)p^2
-
(p^2 + a_p + 1 - p^2)p\\
&&+ (2p^2 + p)+6p(p+1-c_p) + 3(p(p+1-c_p)+c_p) \\
&&+(p^2+2p+1) -2c_p+4p(p+1-b_p + p + 1)\\ 
&=&
p^3 + 20p^2 - (a_pc_p + 9pc_p + 4pb_b) + 20p + 1.
\end{eqnarray*}
which, as predicted, is
the same as the number of points given by (\ref{eqn:predicted_point_count}).

\subsection*{Acknowledgements}

We are grateful to the DFG for
grant Hu 337/5-2 (Schwerpunktprogramm ``Globale Methoden in der
komplexen Geometrie'') .
The second author was partially supported by
the Louisiana Board of Regents Research Competitiveness Subprogram,
contract number LEQSF-(2004-7) RD-A-16.

We would like to thank B.~van Geemen for explaining to us results about
the field of definition of blow ups of varieties defined over $\QQ$,
and C.~Schoen for discussions concerning the computation of $h^{12}$.
Particular thanks also to N.~Yui for numerous discussions on this and related subjects.

\noindent
Klaus Hulek,\\
Institut f\"ur Mathematik (C),\\
Universit\"at Hannover\\
Welfengarten 1, 30060 Hannover, Germany\\
{\tt hulek@math.uni-hannover.de}

\medskip

\noindent
Helena A. Verrill,\\
Department of Mathematics,\\
Louisiana State University\\
Baton Rouge, LA 70803-4918, USA
{\tt verrill@math.lsu.edu}

\end{document}